TWO THEOREMS ON THE STRUCTURE OF
PYTHAGOREAN TRIPLES AND SOME DIOPHANTINE CONSEQUENCES
*By* Konstantine 'Hermes' Zelator

1. **INTRODUCTION**

The major theorem of this paper (Theorem 1) deals with the structure of primitive pythagorean triangles. In Theorem 1, it is proven that there are no primitive pythagorean triangles of the form ($s_1 x^2, s_2 y^2, z$), i.e., triangles with one leg equal to $s_1$ times a perfect square and the other equal to $s_2$ times an integer square, if the integers $s_1$ and $s_2$ satisfy the following conditions:

a) $s_1$ and $s_2$ are both positive odd squarefree integers.
b) $s_1 s_2 = p_1 ... p_n, n \geq 2, p_1 \equiv 5, p_2 \equiv ... \equiv p_n \equiv 1 \mod 8$, where $p_1,..., p_n$ are primes.
c) If $s_1 \equiv 1 \mod 8$, (and then of course $s_2 \equiv 5 \mod 8$), then $s_1$ is a quadratic nonresidue of every divisor $d$ of $s_2$, with $d > 1$.

(If $s_2 \equiv 1 \mod 8$, then the same condition is assumed on $s_2$ with respect to the divisors $d$ of $s_1, d > 1$.)

One can easily find many examples of integers $s_1$, $s_2$ satisfying the above conditions. Below, we offer a few.

1. $s_1$, $s_2$ primes with $s_1 \equiv 1, s_2 \equiv 5 \mod 8$; $s_1$ a quadratic nonresidue of $s_2$.
2. $s_1 = p_1 p_2, p_1 \equiv p_2 \equiv 1 \mod 8, p_1, p_2$ primes, $s_2$ also a prime, $p_1$ a quadratic nonresidue of $s_2$, while $p_1$ a quadratic residue of $s_2$.
3. $s_1 \equiv 1 \mod 8, s_1$ a prime, $s_2 = p_1 p_2, p_1, p_2$ prime with $p_1 \equiv 1, p_2 \equiv 5 \mod 8; s_1$ a quadratic nonresidue of both $p_1$ and $p_2$.
4. $s_1 = p_1 ... p_n; n \geq 3$, $n$ an odd integer, $s_2$ a prime, $s_2 \equiv 5 \mod 8$; $p_1,..., p_n$ primes with $p_1 \equiv ... \equiv p_n \equiv 1 \mod 8$; an even number of integers among the primes $p_1,..., p_n$ being quadratic residues of $s_2$, while the remaining primes (odd in number) being quadratic nonresudues of $s_2$.

Theorem 2 is somewhat more technical in nature. It is shown that if $n$ is a positive odd squarefree integer such that whenever $s_1 s_2 = n$, then $s_1 \equiv 1, s_2 \equiv 5 \mod 8$ (or vice versa), the diophantine equation $(nx^2)^2 + y^4 = z^2$, with $(nx, y) = 1$ is solvable in $Z^+$, if and only if the diophantine equation $d_1 z^2 = d_2^2 x^4 + d_3^2 y^4$ is solvable in $Z^+$, for some positive divisors $d_1, d_2, d_3$ of $n$ with $d_1 d_2 d_3 = n, d_2 < n$ and $d_3 < n$. It is fairly easy to see that $n$, satisfying the hypothesis of Theorem 2, is either the product of an odd number of primes $p \equiv 5 \mod 8$, or otherwise the product of an odd of primes $p \equiv 5 \mod 8$, and any number of primes $q \equiv 1 \mod 8$.

Theorem 3 is simple in its statement and proof. It states that if $p \equiv q \equiv 1 \mod 4$ and $q$ primes, $p$ a quadratic nonresidue of $q$, then the diophantine equation $(px^2)^2 + (qy^2)^2 = z^2$ with $(px, qy) = 1$, has no solution in $Z^+$.

By combining Theorems 2 and 3, Theorem 4 is established. It states that if $p \equiv 1, q \equiv 5 \mod 8$, $p$ a quadratic nonresidue of $q$, the diophantine equation $(pqx^2)^2 + y^4 = z^2$, with $(pqx, y) = 1$, has a solution in $Z^+$, if and only if the equation $pqz^2 = x^4 + y^4$, with $(x, y) = 1$, has a solution in $Z^+$.

2. THEOREM 1

Let $s_1, s_2$ be two odd squarefree positive integers greater than 1, with $(s_1, s_2) = 1$ and such that $s_1 s_2$ is the product of a prime congruent to 5 mod 8 and primes congruent to 1 mod 8, i.e., $s_1 s_2 = p_1 \ldots p_n, n \geq 2, p_1 \equiv 5 \mod 8, p_2 \equiv \ldots \equiv p_n \equiv 1 \mod 8$. In addition, assume that if $s_1 \equiv 1 \mod 8$ (same condition on $s_2$, if $s_2 \equiv 1 \mod 8$), then $s_1$ is a quadratic nonresidue of every divisor d of $s_2$, including $s_2$ itself, with $d > 1$.

Under the above assumptions, the diophantine equation $s_1^2 x^4 + s_2^2 y^4 = z^2$ has no solution in the set of positive integers $Z^+$, with $(s_1 x, s_2 y) = 1$.

*Remark*

One can easily find examples of $s_1, s_2$ satisfying the hypothesis of the theorem. For instance, $s_1, s_2$ being primes with $s_1 \equiv 1, s_2 \equiv 5 \mod 8$ and $s_1$ a quadratic nonresidue of $s_2$. Also, $s_1 = p_1 p_2, p_1 \equiv p_2 \equiv 1 \mod 8$, where $p_1, p_2$ are primes, $s_2$ also a prime, $s_2 \equiv 5 \mod 8$ such that $p_1$ being a quadratic nonresidue of $s_2$, while $p_1$ a quadratic residue of $s_2$; $s_1 \equiv 1 \mod 8$, $s_1$ a prime, while $s_2 = p_1 p_2, p_1, p_2$ primes, $p_1 \equiv 1, p_2 \equiv 5 \mod 8$, and with $s_1$ a quadratic nonresidue of both $p_1$ and $p_2$.

*Proof.*

Let us assume that $(x, y, z)$ is a solution with $(s_1 x, s_2 y) = 1$ in positive integers $x, y, z$ to the equation.

$$(s_1 x^2)^2 + (s_2 y^2)^2 = z^2 \tag{1}$$

In virtue of $(s_1 x, s_2 y) = 1$, we see that equation (1) describes a primitive Pythagorean triangle.

Assume first $x$ to be odd any $y$ even. From (1), it follows that

$$s_1 x^2 = m^2 - n^2, \; s_2 y^2 = 2mn, \; z = m^2 + n^2 \tag{2}$$

For positive integers $m, n$ with $(m, n) = 1$ and $m + n \equiv 1 \mod 2$

The first equation in (2) implies, since $(m, n) = 1$,

$$x = km_1n_1, \quad n = k\frac{|d_1m_1^2 - d_2n_1^2|}{2}, \quad m = k\frac{|d_1m_1^2 + d_2n_1^2|}{2}, \qquad (3)$$

for positive integers $m_1, n_1, d_1, d_2, k$ with $(m_1, n_1) = 1, d_1d_2 = s_1$ and $k = 1$ or 2. (Refer to References [1] or [2].) Note that $s_1$ is squarefree, so must be $d_1$ and $d_2$. Now, $x$ is odd and so, from the first equation in (3), we see that $k = 1$ and $m_1 \equiv n_1 \equiv 1 \bmod 2$.

On the other hand, the second equation in (2) implies the following two possibilities or sub cases:

$$m = 2\delta_1 M^2, n = \delta_2 N^2 \qquad (4)$$

$$m = \delta_1 M^2, n = 2\delta_2 N^2 \qquad (5)$$

Where $\delta_1, \delta_2$ are positive integers with $\delta_1\delta_2 = s_2$, (they are squarefree since $s_2$ is), and $M, N$ are positive integers with $(M, N = 1)$, since $(m, n) = 1$.

Also, since, as we have shown above, k=1, the last two equations in (3) take the form

$$2n = |d_1m_1^2 - d_2n_1^2|, \quad 2m = d_1m_1^2 + d_2n_1^2 \qquad (6)$$

First assume (4) to be the case. Then of course, since m and n have different parities, m must be even and n odd, and so N must also be odd.

Combining equations (4) and (6), we obtain

$$2\delta_2 N^2 = |d_1m_1^2 - d_2m_1^2|, \quad 4\delta_1 M^2 = d_1m_1^2 + d_2n_1^2. \qquad (7)$$

By virtue of $N \equiv m_1 \equiv n_1 \equiv 1 \bmod 2$ (7) yields

$$2\delta_2 \equiv |d_1 - d_2| \bmod 8, \quad 4\delta_1 M^2 \equiv d_1 + d_2 \bmod 8. \qquad (8)$$

According to the hypothesis of Theorem 1, either $s_1$ is a prime congruent to $5 \bmod 8$ or a product of a prime congruent to $5 \bmod 8$ and primes congruent to $1 \bmod 8$ and $s_2$ is a product of primes congruent to $1 \bmod 8$, or vice versa. Suppose, under the hypothesis of Theorem 1, that $s_2 \equiv 5 \bmod 8$, while $s_1 \equiv 1 \bmod 8$; then from $d_1d_2 = s_1$ and $\delta_1\delta_2 = s_2$, we conclude $d_1 \equiv d_2 \equiv 1 \bmod 8$ and $\delta_1 \equiv 1, \delta_2 \equiv 5$ or $\delta_1 \equiv 5, \delta_1 \equiv 1 \bmod 8$. Then, however, the first congruence in (8) implies

$2\delta_2 \equiv 0 \bmod 8$, a contradiction since $\delta_2$ is odd.

Now, if $s_1 \equiv 5$ and $s_2 \equiv 1 \bmod 8$, we obtain, from $s_1 = d_1d_2$ and $s_2 = \delta_1\delta_2$, under the hypothesis of the theorem,

$d_1 \equiv 1, d_2 \equiv 5$ or $d_1 \equiv 5, d_2 \equiv 1$ and $\delta_1 \equiv \delta_2 \equiv 1 \mod 8$

The first equation in (8) implies $2\delta_2 \equiv 4 \mod 8$, a contradiction since $\delta_2$ is odd.

Next, assume equation (5) to be the case; then M is odd. By combining equations (5) and (6), we arrive at

$$4\delta_2 N^2 = \pm(d_1 m_1^2 - d_2 n_1^2), \quad 2\delta_1 M^2 = d_1 m_1^2 + d_2 n_1^2. \tag{9}$$

Repeating the reasoning we applied when we considered (8), we see that if $s_1 \equiv 5$ and $s_2 \equiv 1 \mod 8$, then $d_1 \equiv 1, d_2 \equiv 5$ or $d_1 \equiv 5, d_2 \equiv 1 \mod 8$ while $\delta_1 \equiv \delta_2 \equiv 1 \mod 8$.

The second of (9) yields $2\delta_1 M^2 \equiv 2\delta_1 \equiv d_1 + d_2 \mod 8$; but $2\delta_1 \equiv 2$, while $d_1 + d_2 \equiv 6 \mod 8$, whence a contradiction. Finally, assume $s_1 \equiv 1$ and $s_2 \equiv 5 \mod 8$.

According to Legendre's theorem, since $Nm_1 n_1 \neq 0$, the first equation in (9) implies that $d_1 \cdot d_2$ is a quadratic residue of $\delta_2$; but $d_1 d_2 = s_1 \equiv 1 \mod 8$ and $\delta_2$ divides $s_2$. Therefore, if $\delta_2$ is a divisor of $s_2$, with $\delta_2 > 1$, by the hypothesis of the theorem, it would follow that $s_1$ is a quadratic nonresidue of $\delta_2$, whence a contradiction.

Now, if $\delta_2 = 1$ then $\delta_1 \delta_2 = s_2 \Rightarrow \delta_1 = s_2$.

Then, however, the second equation of (9) and Legendre's theorem imply that $-d_1 d_2 = -s_1$ is a quadratic residue of $\delta_1 = s_2$; but -1 is a quadratic residue of $s_2$ (since $s_2$ is a product of primes congruent to 1 mod 4). Hence $s_1$ is a quadratic residue of $s_2$, contrary to the hypothesis of the theorem.

To conclude the proof of Theorem 1, let us go back to (1) and consider the case where x is even and y odd. Unless we use the special assumption that if $s_1 \equiv 1 \mod 8$, then $s_1$ is a quadratic nonresidue of every divisor if $s_2$ (except 1), it is clear that the case x = even and y = odd is identical in treatment with the case x = odd and y = even which we have already treated.

Now, if $s_1 \equiv 1 \mod 8$, then $s_2 \equiv 5 \mod 8$; then the treatment of the case x = even, y = odd, $s_1 \equiv 1$ and $s_2 \equiv 5 \mod 8$, is obviously identical with the treatment of the case x = odd, y = even, $s_1 \equiv 5$ and $s_2 \equiv 1 \mod 8$, which has already been done.

## 3. THEOREM 2

Let *n* be a squarefree positive integer such that if $s_1, s_2$ are any divisors of *n* with $s_1 s_2 = n$, then $s_1 \equiv 1, s_2 \equiv 5 \mod 8$ or vice versa. Then the diophantine equation $n^2 x^4 + y^4 = z^2$ with $(nx, y) = 1$, has a solution in the set of positive integers $Z^+$, if and only if, the diophantine equation $d_1 z^2 = d_2^2 x^4 + d_3^2 y^4$ has a solution, with $(d_2 x, d_3 y) = 1$, in positive integers $x$, $y$, $z$, for some positive divisors $d_1, d_2, d_3$ of $n$, such that $d_1 d_2 d_3 = n$, with $d_2 < n$ and $d_3 < n$.

*Remark.*

It follows from the hypothesis of the theorem that $n \equiv 5 \mod 8$ (take $s_1 = n$ and $s_2 = 1$); one easily finds examples of numbers n that satisfy the hypothesis.

1. $n =$ prime, $n \equiv 5 \mod 8$.
2. $n = p_1 p_2$, $p_1, p_2$ primes with $p_1 \equiv 5, p_1 \equiv 1 \mod 8$.
3. $n = p_1 p_2 p_3 ... p_i$ $n \geq 3$; $p_1 \equiv 5, p_2 \equiv p_3 \equiv ... \equiv p_i \equiv 1 \mod 8$.
4. *n* is the product of an odd number of primes congruent to 5 mod 8 and any number of primes congruent to 1 mod 8.
5. In fact, *n* is either the product of an odd number of primes $p \equiv 5 \mod 8$, or otherwise the product of an odd number of primes $p \equiv 5 \mod 8$ and any number of primes $q \equiv 1 \mod 8$.

*Proof.*
Let $x$, $y$, $z$ be positive integers, with $(nx, y) = 1$, satisfying the equation

$$(nx^2)^2 + (y^2)^2 = z^2 \tag{10}$$

We distinguish between the cases ($x =$ even, $y =$ odd) and ($x =$ odd, $y =$ even).

*Case 1;* $x =$ odd, $y =$ even

Equation (10), together with the condition $(nx, y) = 1$, imply

$$nx^2 = r^2 - t^2, \quad y^2 = 2rt, \quad x = r^2 + t^2 \tag{11}$$

For positive integers $r, t$ with $(r, t) = 1$ and $r + t \equiv 1 \mod 2$.

The second equation (11) implies

$$\left. \begin{array}{l} r = R^2, \; t = 2T^2 \\ \text{or} \\ r = 2R^2, \; t = T^2 \end{array} \right\} \quad (R, T) = 1$$

The latter possibility is ruled out, for if it holds, the first equation in (11) gives

$$nx^2 \equiv 4R^4 - T^4 \mod 8 \tag{12}$$

But T is odd (since $t = T^2, r \equiv 0 \mod 2$ and $r + t \equiv 1 \mod 2$), hence (12) implies $nx^2 \equiv 4R^4 - 1 \equiv 3 \text{ or } 7 \mod 8$; however, x is odd and so $nx^2 \equiv n \equiv 3 \text{ or } 7 \mod 8$, contrary to the fact that $n \equiv 5 \mod 8$ (refer to the remark underneath Theorem 2).

Now, suppose that $r = R^2, t = 2T^2$. The first equation in (11) implies, from References [2] or [1], (and since $(r,t)=1$),

$$x = k \cdot R_1 T_1, \quad t = \frac{k \cdot |s_1 R_1^2 - s_2 T_1^2|}{2}, \quad r = \frac{k(s_1 R_1^2 + s_2 T_1^2)}{2} \tag{13}$$

for integers $R_1, T_1$ with $(R_1, T_1) = 1$ and integers $s_1, s_2$ with $s_1 s_2 = n$ and also k = 1 or 2.

Now, since $x$ is odd it follows from the first equation of (13), $k=1$ and $R_1 \equiv T_1 \equiv 1 \mod 2$.
Since $r = R^2$, the third equation in (13) implies

$$2R^2 = s_1 R_1^2 + s_2 T_1^2,$$
$$2R^2 \equiv s_1 R_1^2 + s_2 T_1^2 \equiv s_1 + s_2 \mod 8.$$

However, R is odd (recall we are in the subcase $r = R^2, t = 2T^2$), and so the last congruence gives $s_1 + s_2 \equiv 2 \mod 8$; however, we have $s_1 s_2 = n$ and so, according to the hypothesis of the theorem, we must have $s_1 \equiv 1$, $s_2 \equiv 5 \mod 8$ or vice versa. At any rate, we obtain $s_1 + s_2 \equiv 6 \mod 8$, contradicting the congruence $s_1 + s_2 \equiv 2 \mod 8$ obtained above.

Note that since $s_1 - s_2 \equiv 4 \mod 8$ and $t = 2T^2$ and therefore the second equation in (13) cannot be rendered impossible via a congruence modulo 8.

*Case 2:* $x$ = even, $y$ = odd.

Equation (1) gives

$$nx^2 = 2rt, \quad y^2 = r^2 - t^2, \quad z = r^2 + t^2 \tag{14}$$

for positive integers r, t with $(r,t)=1$ and $r+t \equiv 1 \mod 2$

The second equation of (14) shows, since y is odd, that r is odd and t even (consider it mod 4). Thus, the first equation in (14) implies (in virtue of $(r,t)=1$)

$$r = d_1 R^2, \quad t = 2 d_2 T^2 \tag{15}$$
with $d_1 d_2 = n (d_1, d_2 > 0)$.

On the other hand, from $y^2 = r^2 - t^2$ and $(r,t) = 1$, we obtain in
$$r \equiv R_1^2 + T_1^2, \quad t = 2 R_1 T_1 \tag{16}$$

for positive integers $R_1, T_1$ with $(R_1, T_1) = 1$ and $R_1 + T_1 \equiv 1 \mod 2$. The second equations of (16) and (15) give

$$2d_2T^2 = 2R_1T_1 \ , \ (R_1,T_1)=1 \tag{17}$$

Therefore
$$R_1 = d_3R_2^2 \ , \ T_1 = d_4T_2^2 \tag{18}$$
and with $d_3d_4 = d_2$, $(d_3,d_4 > 0)$.

By using the first of (15), (16) and equations (18), we obtain

$$d_1R^2 = d_3^2R_2^4 + d_4^2T_2^4 \tag{19}$$

with $(d_3R_2, d_4T_2)=1$ (since $(R_1,T_1)=1$). Note that $d_1d_3d_4 = n$, since $d_1d_2 = n$ and $d_3d_4 = d_2$. We claim that $d_3 < n$ and $d_4 < n$ in (19). For if, say $d_3 = n$, then from $d_1d_3d_4 = n$, it follows that $d_1 = d_4 = 1$ (remember $d_1,d_2,d_3 > 0$). Also, since $d_2 = d_3d_4$, we must have $d_2 = n$.

From $d_3 = n$, $d_4 = 1$ and equation (18), we obtain

$$T_1 = T_2^2 \ , \ R_1 = nR_2^2 \tag{20}$$

From $d_2 = n$ and equation (17), we obtain

$$T^2 = R_2^2T_2^2 \tag{21}$$

From the second equation of (16), we arrive at

$$t = 2nR_2^2T_2^2 \tag{22}$$

Combining (21) and (22), we obtain

$$t = 2nT^2 \tag{23}$$

On the other hand, $d_1 = 1$, $d_2 = n$ (15) and (16) imply $x^2 = R^2T^2$ and so by (21) we obtain

$$x^2 = R^2R_2^2T_2^2 \tag{24}$$

Equation (24) shows that $s^2 \geq T_2^2$. We claim that the equal sign cannot hold, for if $x^2 = T_2^2$, then (24) implies $R^2 = R_2^2 = 1$. However, (19) is then rendered impossible for it gives, on account of $d_1 = 1$, $d_3 = n$, $d_4 = 1$, $1 = n^2 + T_2^4$, which is of course impossible.

Hence (24) shows that it mush be $0 < T_2^2 < x^2$, and so we are led to an indefinite descent with respect to the initial equation (1). The same argument is made if we assume $d_4 = n$ in (19).

Finally, to construct a solution to (1), from a solution to (19) is an easy matter. Given a solution $(R, R_2, T_2)$ in $Z^+$ to (19), with $d_1 d_3 d_4 = n$, $d_3 < n$, $d_4 < n$, and $(d_3 R_2, d_4 T_2) = 1$, one can construct a solution $(x, y, z)$ in $Z^+$ to (19), by simply tracing back through equations (19) to (14)

## 4. THEOREM 3

Let $p$, $q$ be primes with $p \equiv q \equiv 1 \bmod 4$. Assume that $p$ is a quadratic nonresidue of $q$ (and thus q is a nonresidue of p, by the reciprocity law). Then the diophantine equation $(px^2)^2 + (qy^2)^2 = z^2$, with $(px, qy) = 1$ has no solution in the set of positive integers $Z^+$.

*Proof.* Assume $x, y, z$ in $Z^+$, to satisfy

$$(px^2)^2 + (qy^2)^2 = z^2 \tag{25}$$

With $(px, qy) = 1$ Then equation (25) describes a primitive pythagorean triangle, and so by assuming x to be odd and y even (without any loss of generality), we obtain

$$px^2 = m^2 - n^2, \quad qy^2 = 2mn, \quad z = m^2 + n^2 \tag{26}$$

for positive integers $m$, $n$ with $(m, n) = 1$ and $m + n \equiv 1 \bmod 2$. Since $x$ is odd and $p \equiv 1 \bmod 4$, we have $px^2 \equiv 1 \bmod 4$, and so a congruence modulo 4 in the first of (26) show that $m \equiv 1 \bmod 2$ and $n \equiv 0 \bmod 2$.

Then the second equation of (26) yields two possibilities, since $(m, n) = 1$, either

$$\text{or} \quad \left.\begin{array}{c} m = qm_1^2, \ n = 2n_1^2 \\ \\ m = m_1^2, \ n = 2qn_1^2 \end{array}\right\} \quad (m, n) = 1 \tag{27}$$

If the first possibility of (27) holds, then the first equation in (26) implies

$$px^2 \equiv -n^2 \bmod q \tag{28}$$

But, -1 is a quadratic residue of $q$ (since $q \equiv 1 \bmod 4$), and so by virtue of $n \not\equiv 0 \bmod q$ (because $m \equiv 0 \bmod q$ and $(m, n) = 1$, it follows from (27) that $p$ is a quadratic residue of $q$, contradicting the hypothesis of the theorem.

If the second possibility of (27) is the case, we then obtain from the first of (26), $px^2 \equiv m^2 \bmod q$, and a similar contradiction, as previously, is obtained.

## 5. THEOREM 4

Let $p, q$ be primes with $p \equiv 1$ and $q \equiv 5 \mod 8$. Also, suppose that $p$ is a quadratic nonresidue of $q$, (and so $q$ is a quadratic nonresidue of $p$).

Under the above assumption, the diophantine equation $(pqx^2)^2 + y^4 = z^2$, with $(pqx, y) = 1$, has a solution in positive integers $x, y, z$ if and only if the diophanitne equation $pqz^2 = x^4 + y^4$ has solution, with $(x, y) = 1$, in positive integers $x, y, z$.

*Proof.*

Clearly, the hupothesis of this theorem satisfies the hypothesis of Theorem 2, with $n = pq$.

Hence, according to Theorem 2, $(pqx^2)^2 + y^4 = z^2$ has a solution, with $(pqx, y) = 1$, in $\mathbb{Z}^+$, if and only if the equation

$$d_1 z^2 = d_2^2 x^4 + d_3^2 y^4 \tag{29}$$

has a solution in $\mathbb{Z}^+$, for some $d_1, d_2, d_3$ with $(d_2 x, d_3 y) = 1$, $d_1 d_2 d_3 = n$, $d_2 < n$ and $d_3 < n, (d_1, d_2, d_3 > 0)$.

But $n = pq$, $pq = d_1 d_2 d_3$, thus if $d_1 = 1$, then on account of $d_2, d_3 < pq$, it must be $d_2 = p$, $d_3 = q$ or vice versa of course.

However, this would imply, by equation (29), an equation
$z^2 = p^2 x^4 + q^2 y^4$ with $(px, qy) = 1$, Which is, by Theorem 3, impossible in $\mathbb{Z}^+$.

Hence we see that we cannot have $d_1 = 1$ in (29). Now, if $d_1 = p$, then $d_2 = q$ and $d_3 = 1$ or $d_2 = 1$ and $d_3 = q$. So (29) would yield

$$pz^2 = q^2 x^4 + y^4, \quad (qx, y) = 1$$

or $\tag{30}$

$$pz^2 = x^4 + q^2 y^4, \quad (x, qy) = 1$$

However, under the conditions $(qx, y) = 1, (x, qy) = 1$, the equations in (30) imply that p is a quadratic residue of q, contradicting the hypothesis of the theorem.

A similar argument is left for the case of $d_1 = q$. Consequently, we see that (29) may have a solution $\mathbb{Z}^+$, only for $d_1 = pq$ and $d_2 = d_3 = 1$.

Thus, the equation $(pqx^2)^2 + y^4 = z^2$, with $(pqx, y) = 1$, is equivalent in $\mathbb{Z}^+$, with the equation $pqz^2 = x^4 + y^4$, $(x, y) = 1$.